\documentclass[10pt,a4paper]{article}
\usepackage{graphicx}
\usepackage{amsmath}
\usepackage{amsfonts}
\usepackage{amssymb}
\usepackage{amsthm}
\begin{document}

\def\R{\mathbb{R}}
\def\H{\mathcal{H}}
\def\d{\textrm{div}}
\def\v{\textbf{v}}
\def\w{\textbf{w}}
\def\u{\textbf{u}}

\newtheorem{lem}{Lemma}
\newtheorem{prop}{Proposition}
\newtheorem{cor}{Corollary}
\newtheorem{thm}{Theorem}
\newtheorem{ass}{Assumption}
\newtheorem{definition}{Definition}
\newtheorem{con}{Conjecture}
\newtheorem{claim}{Claim}
\newtheorem{remark}{Remark}

\title{Almost Everywhere Regularity for the Free Boundary of the Normalized p-harmonic Obstacle problem $p>2$.}
\author{John Andersson}

\maketitle
\begin{abstract}
Let $u$ be a solution to the normalized p-harmonic obstacle problem with $p>2$. That is, $u\in W^{1,p}(B_1(0))$, $2<p<\infty$, $u\ge 0$ and
$$
\d\left( |\nabla u|^{p-2}\nabla u\right)=\chi_{\{u>0\}}\textrm{ in }B_1(0)
$$
where $u(x)\ge 0$ and $\chi_A$ is the characteristic function of the set $A$.
Our main result is that for almost every free boundary point, with respect to the $(n-1)-$Hausdorff measure, there is a neighborhood where the free boundary
is a $C^{1,\beta}-$graph. That is, for $\H^{n-1}-$a.e. point  $x^0\in \partial \{u>0\}\cap B_1(0)$ there is an $r>0$
such that $B_r(x^0)\cap \partial \{u>0\}\in C^{1,\beta}$.
\end{abstract}

\section{Introduction.}
In this article we will consider the following $p-$harmonic obstacle problem
\begin{equation}\label{main}
\begin{array}{ll}
 \d\left( |\nabla u|^{p-2}\nabla u\right)=\chi_{\{u>0\}} &  \textrm{ in }D \\
 u(x)=g(x)\ge 0 & \textrm{ on } \partial D
\end{array}
\end{equation}
where $g(x)\ge 0$ is the restriction of a $W^{1,p}(D)$ function to $\partial D$. We will denote $\Omega_u=\{u>0\}$ and 
the free boundary $\Gamma =\partial \Omega \cap D$.

A solution to (\ref{main}) can be found by minimizing 
$$
\int_D \frac{1}{p}|\nabla u|^p + \max(u,0) dx
$$
in $K=\{u\in W^{1,p}(D);\; u-g\in W^{1,p}_0(D)\}$. It is well known that solutions to (\ref{main}) are $C^{1,\alpha}_{\textrm{loc}}(D)$ for some
$\alpha>0$, see for instance \cite{CL}.

Our main interest will be the free boundary $\Gamma =\partial \Omega \cap D$. In the special, linear, case $p=2$ it is known that 
the free boundary is a $C^{1,\beta}$ (even analytic) graph around almost every point, with respect to $\H^{n-1}$, of the free boundary
\cite{CaffRevisited}. However, the techniques to prove regularity of the free boundary in the linear case breaks down when $p\ne 2$. 
In particular, it seems difficult to use the strong comparison and boundary comparison methods of \cite{CaffRevisited} for 
$p-$harmonic problems. Moreover, due to the nonlinearity of the $p-$laplace operator, the solution, $u$, solves a different PDE compared to the 
derivative of the solution $\frac{\partial u}{\partial x_i}$. New techniques are needed to investigate the free boundary for the $p-$harmonic obstacle problem for $p\ne 2$.
In this article we develop these techniques for $p>2$.

We will use a flatness improvement argument and assume that the solution is almost one dimensional:
\begin{equation}\label{BigFlatAss}
 \|\nabla' u\|_{L^2_A(B_1)}\le \epsilon_0
\end{equation}
where $L^2_A$ is the normal $L^2$ space with a suitable weight and \\ $\nabla'=(\partial_1,\partial_2,...,\partial_{n-1},0)$.
We prove that when (\ref{BigFlatAss}) is satisfied for an $\epsilon_0$ small enough then, for some $\tau<1$,
\begin{equation}\label{BetterFlatness}
 \left\|\nabla' \frac{u(sx)}{s^{p/(p-1)}}\right\|_{L^2_A(B_1)}\le \tau\|\nabla' u\|_{L^2_A(B_1)}.
\end{equation}
A standard iteration of (\ref{BetterFlatness}) implies that the free boundary is $C^{1,\beta}$ in $B_{1/2}$. 

The ideas we use to prove the $C^{1,\beta}$ regularity of the free boundary are standard flatness improvement arguments
from geometric measure theory and non-linear systems in the calculus of variations. However, the techniques will look somewhat 
different since we are dealing with regularity of the free boundary and not with the regularity of a function.

We will assume that we have a sequence of solutions $u^j$ to (\ref{main}) in $B_1(0)$
satisfying (\ref{BigFlatAss}) with constants $\epsilon_j\to 0$. We then show, more or less, that the limit function $u^0=\lim_{j\to \infty}u^j$
satisfies an appropriate estimate of the kind (\ref{BetterFlatness}). It is then a matter of showing that the limit $u^j\to u^0$
is strong enough to draw the conclusion that (\ref{BetterFlatness}) holds for all $u^j$ with small enough $\epsilon_j$.
For details and precise statements we refer the reader to the main text.

The second main result of this paper is that (\ref{BigFlatAss}) is actually satisfied at almost every free boundary point, $x^0\in \Gamma$, in a 
small enough ball $B_{r(x^0)}(x^0)$. This result is more subtle 
than it appears to be. We know, \cite{LS} ($p>2$, see also \cite{Rodr} for these results in a much more general situation), that $\Gamma$
has finite $(n-1)-$Hausdorff measure so a famous result of de Giorgi implies that the free boundary has a measure theoretic normal
at a.e. free boundary point. This means that $u^0=\lim_{r\to 0}\frac{u(rx+x^0)}{r^{p/(p-1)}}$ will converge to a solution to (\ref{main})
in $\R^n$ with $\overline{\{u^0>0\}}=\{x; x\cdot \eta \ge 0\}$, that is $\Omega_{u^0}$ is a half-space and the free boundary a hyperplane.
However, since the equation is not uniformly elliptic we can not use the Kauchy-Kovalevskaya Theorem and conclude that $u^0$ is one
dimensional. There should be such a simple argument to show that $u^0$is one dimensional - but we could not find any. Instead we will construct a 
Carleson measure related to 
$u^0$ and use that to show that there is a blow-up of $u^0$ that is one dimensional. The argument is interesting in its own right and based 
on an idea from \cite{DavSem}.
A Lemma from geometric measure theory (see for instance \cite{DeLell}), used in the 
free boundary context in \cite{Bellman}, then implies that a blow-up of $u^0$ is in fact also a blow-up of $u$ at $x^0$ - at least for 
almost every $x^0$ in the free boundary.

Our main Theorem is

\vspace{2mm}

{\bf Theorem \ref{MainResult}.} 
{\sl Let $u$ be a solution to the $p-$harmonic obstacle problem, $p\ge 2$, in $B_1(0)$ then there exists an open set $\Gamma_0\subset \Gamma$
 such that $\H^{n-1}(\Gamma\setminus \Gamma_0)=0$ and for every $x^0\in \Gamma_0$ there exists an $r=r(x^0)$ such that 
 $\Gamma_u\cap B_r(x^0)$ is a $C^{1,\alpha}$ graph.}

\vspace{2mm}

The plan of the paper is as follows. In the next section we gather some known regularity results for $p-$harmonic
obstacle problem. In section \ref{sec:BlowUp} we recall some results relation to blow-ups of solutions. In section
\ref{sec:Carleson} we construct the Carleson measure that is vital to show that (\ref{BigFlatAss})
is satisfied in a small ball centered at a.e. free boundary point. In the following section we show that any global
solution with the free boundary being a hyperplane has a blow-up that is one dimensional. In section \ref{sec:lin}
we show that we can linearize the problem and that the convergence of the linearizing sequence is strong enough. We are then 
ready to prove Theorem \ref{MainResult}.

\subsection{Notation.} We will use the following notation:
\begin{enumerate}
 \item We will use $B_r(x^0)=\{x\in \R^n;\; |x-x^0|<r\}$, $B_r^+(x^0)=\{x\in \R^n;\; |x-x^0|<r,\; x_n\ge 0\}$.
 \item An $'$ will indicate that the $x_n-$coordinate is excluded: $x'=(x_1,...x_{n-1},0)$, $B_r'(z')=\{x';\; |x'-z'|<r\}$,
 $\nabla'=(\partial_1,..., \partial_{n-1},0)$ et.c.
 \item On occasion we are going to use, for a unit vector $\eta\in \R^n$, $\nabla'_\eta=\nabla -\eta \left(\eta\cdot \nabla\right)$
 (the gradient operator restricted to the subspace orthogonal to $\eta$).
 \item $W^{1,p}(D)$ will be the usual Sobolev space and $W^{1,p}_A(D)$ will denote the Sobolev space with weight $A(x)$; That is the space of all 
 weakly differentiable functions with norm $\|u\|_{W_A^{1,p}(D)}^p=\int_D A(x)|u|^pdx+\int_D A(x)|\nabla u|^pdx$.
 \item By $\Omega$ and $\Gamma$ we denote the non-coincidence set $\Omega=\{x;\; u(x)>0\}$ and the free boundary $\Gamma=\partial \Omega\cap D$.
 At times we will write $\Gamma_u$ and $\Omega_u$ to indicate the function whose free boundary we are considering.
 \item A blow-up of $u$ at $x^0\in D$ will be the limit of any convergent sequence $\frac{u(r_j x+x^0)}{r_j^{p/(p-1)}}$ as $r_j\to 0$.
 \item Characteristic function of the set $S$ will be denoted by $\chi_S$.
\end{enumerate}

\section{Regularity and known results.}\label{sec:reg}

In this section we gather some known regularity results for p-harmonic obstacle problems.

The following regularity results are well known for solutions to the $p$-harmonic obstacle problem.

\begin{lem}\label{RegularityPharmObst}
 Let $u$ be a solution to the $p-$harmonic obstacle problem in $B_2(0)$ then there exists an $\alpha>0$ and a constant $C$ such that
 \begin{equation}\label{C1alphaReg}
 \|u\|_{C^{1,\alpha}(B_{3/2}(0))}\le C\|u\|_{L^\infty(B_2(0))}.
 \end{equation}
 Furthermore there exists constants $C,c>0$ such that if $x^0\in \Gamma\cap B_1(0)$ then
 \begin{equation}\label{OptGrowth}
 cr^{p/(p-1)}\le \sup_{B_r(x^0)}u\le Cr^{p/(p-1)}.
 \end{equation}
\end{lem}
The $C^{1,\alpha}-$regularity can be found in \cite{Lewis}. The growth estimate (\ref{OptGrowth}) was proved in \cite{KKPS} or \cite{Rodr}.


We will also need a weak regularity result which can be proved by a standard difference quotient argument. Later we will use the same ideas to 
prove Lemma \ref{w22estimates} and since the proofs are very similar we do not prove the Lemma here.

\begin{lem}\label{L2}
Assume that $u$ is a solution to the p-harmonic obstacle problem in $\tilde{D}$ and $D\subset\subset \tilde{D}$. Let
$\delta=\textrm{dist}(D,\tilde{D}^c)$ and $\psi$ be a standard cut off function: $\psi\in C^\infty_c(\tilde{D})$, $\psi=1$ in $D$
and $|\nabla \psi|\le \frac{C_0}{\delta}$. 

Then there exists a constant $C_1=C_1(n,p)$ depending 
only on the dimension and $p$ such that 
$$
\int_{\tilde{D}}\psi^2 |\nabla u|^{p-2}\left| \frac{\partial \nabla u}{\partial x_i}\right|^2dx \le \frac{C_1 C_0^2}{\delta^2}\int_{\tilde{D}\setminus D}|\nabla u|^{p-2}\left|\frac{\partial u}{\partial x_i}\right|^2dx.
$$
\end{lem}

\subsection{Eigenfunction Expansions for the Linearized problem.}

Consider the Rayleigh functional  
\begin{equation}\label{Rayleigh}
J(u)=\frac{\int_{\partial B_1^+}|x_n|^{\frac{p-2}{p-1}}|\nabla_\phi u|^2}{\int_{\partial B_1^+}|x_n|^{\frac{p-2}{p-1}}u^2}
\end{equation}
defined on all functions in $W^{1,2}(\partial B_1^+, |x_n|^{(p-2)/(p-1)})$ that vanish on $\{x_n=0\}$. Here $\nabla_\phi$
is the gradient restricted to the sphere.

Notice that there is not any problem to define the boundary values of $u\in W^{1,2}(\partial B_1^+, |x_n|^{(p-2)/(p-1)})$
in the sense of traces since since $W^{1,1}(\partial D_1^+(0))$ has a trace operator and for every $u\in W^{1,2}(\partial B_1^+, |x_n|^{(p-2)/(p-1)})$
$$
\int_{\partial B_1^+}|\nabla_\phi u|dx\le 
\left( \int_{\partial B_1^+}|x_n|^{-\frac{p-2}{p-1}}\right)^{\frac{1}{2}}\left( \int_{\partial B_1^+}|x_n|^{\frac{p-2}{p-1}}|\nabla u|^2\right)^{\frac{1}{2}}<
\infty.
$$

A minimizer to (\ref{Rayleigh}) may be extended homogeneously to a global solution in $\R^n_+$ to
$$
\begin{array}{ll}
 \d\left(|x_n|^{\frac{p-2}{p-1}}\nabla u\right)=0 & \textrm{ in }\R^n_+ \\
 u(x)=0 & \textrm{ on } x_n=0.
\end{array}
$$

Using a standard argument from functional analysis, see appendix D in \cite{Evans}, gives the following lemma.

\begin{lem}\label{eigfunclem}
 Let $v$ be a solution to 
 \begin{equation}\label{eigenfunction}
 \begin{array}{ll}
 \d\left( |x_n|^{\frac{p-2}{p-1}}\nabla v\right)=0 & \textrm{ in } B_1^+(0) \\
 v(x)=0 & \textrm{ on } \{x_n=0\}
 \end{array}
 \end{equation}
 then there is a sequence of solutions $q_j$ to (\ref{eigenfunction}) such that
 \begin{enumerate}
  \item $q_j$ is $\lambda_j-$homogeneous,
  \item $\frac{1}{p-1}=\lambda_1<\lambda_2\le \lambda_3\le \lambda_4\le....$,
  \item $\int_{B_1^+}|x_n|^{\frac{p-2}{p-1}}\nabla q_i\cdot \nabla q_j =0$ for $i\ne j$,
  \item $v(x)=\sum_{j=1}^\infty a_jq_j(x)$ for some constants $a_j$.
 \end{enumerate}
\end{lem}

\section{Blow-ups}\label{sec:BlowUp}

In this section we analyze the blow-ups of the solutions. One of the main goals of this paper is to show that blow-ups are unique, that is that $u^0$
is independent of the sequence $r_j\to 0$. In this section we will gather some known results. But we will begin by define what we mean by a blow-up.

\begin{definition}
 If $u$ is a solution to the $p-$harmonic obstacle problem in $B_1(0)$ and $x^0\in \Gamma_u$. Then we say that $u^0(x)$ is a blow-up
 of $u$ at $x^0$ if there is a sequence $r_j\to 0$ such that
 $$
 \lim_{j\to \infty}\frac{u(r_j x+x^0)}{r_j^{\frac{p}{p-1}}}=u^0(x),
 $$
 where the limit is considered weakly in $W^{1,p}_{loc}(\R^n)$.
 
 We denote the set of blow-ups of $u$ at $x^0$ by $Blo(u,x^0)$. That is
 $$
 Blo(u,x^0)=\left\{ u^0(x); \; \exists r_j\to 0,\; u^0(x)=\lim_{j\to \infty}\frac{u(r_j x+x^0)}{r_j^{\frac{p}{p-1}}}\right\}.
 $$
\end{definition}

The following properties (in Lemma \ref{ThisWouldNeedAName} and \ref{ThisToo}) are well known for blow-ups.

\begin{lem}\label{ThisWouldNeedAName}
 Suppose that $u$ solves the $p-$harmonic obstacle problem in $B_1(0)$ and $x^0\in \Gamma_u$. Then for 
 any sequence $r_j\to 0$ there exists a subsequence $r_{j_k}\to 0$ such that
 $$
  \lim_{k\to \infty} \frac{u(r_{j_k} x+x^0)}{r_{j_k}^{\frac{p}{p-1}}}
 $$
 exists. That is $Blo(u,x^0)$ is not empty. If we call the limit $v(x)$ then $v\ne 0$ and $v$ is a solution to the 
 $p-$harmonic obstacle problem.
\end{lem}

Again see \cite{KKPS} or \cite{Rodr} for a proof.

\begin{lem}\label{ThisToo}
 Let $u$ be a solution to (\ref{main}). Then
 \begin{enumerate}
  \item The free boundary $\Gamma_u$ has locally finite $(n-1)-$Hausdorff measure $\H^{n-1}$.
  \item The set $\{u>0\}$ has locally finite perimeter.
  \item For $\H^{n-1}$-a.e. point $x^0\in \Gamma_u$ there exists a measure theoretic normal $\eta_{x^0}$ such that
  $$
  \Omega_{r}=\left\{ \frac{u(rx+x^0)}{r^{\frac{p}{p-1}}}>0\right\}\to \left\{ x,\; x\cdot \eta \ge 0\right\}
  $$
  in the sense of Hausdorff convergence.
 \end{enumerate}
\end{lem}

The first two points can be found in \cite{LS} or in \cite{Rodr} for the last point see \cite{Giusti}.

The following Lemma was proved in \cite{Bellman} for a slightly different problem. But the proof is line for line the same 
for $p$-harmonic obstacle problems, there are some slight typos in the proof in \cite{Bellman} and the reader might 
benefit from considering the forthcoming paper \cite{Linearization} where the some more details are given. 
The proof is based on ideas from geometric measure theory, see for instance \cite{DeLell}.

\begin{lem}\label{BlowOfBlowAreBlow}
 Suppose that $u$ is a solution to the $p-$harmonic obstacle problem. Then for $\H^{n-1}-$almost every free boundary point $x^0$
 it holds that if $u^{x^0}\in Blo(u,x^0)$ and $x^1\in \Gamma_{u^{x^0}}$ then
 $$
 \lim_{r\to 0}\frac{u^{x^0}(rx+x^1)}{r^{\frac{p}{p-1}}}\in Blo(x^0,u).
 $$
\end{lem}

{\bf Remark:} Informally the Lemma states that {\sl ``Blow-ups of blow-ups at a point $x^0$ are also blow-ups at $x^0$.''}

\section{A Carleson measure.}\label{sec:Carleson}

In this section we prove a simple, but technical, lemma. Later it will be used to prove that any global solution,
with free boundary being a hyperplane, has a blow-up that is one dimensional. In particular, if the functions $g^j$
in Lemma (\ref{DiniVan}) are derivatives $\frac{\partial u}{\partial x_j}$, $j=1,2,...,n-1$, then the lemma states that there is 
a blow-up sequence $\frac{u(r_k x+ y^k)}{r_k^{p/(p-1)}}\to u^0$ such that $u^0$ is independent of the $x'$ variables 
(see Proposition \ref{flatnessprop}). 

The proof is very similar to a proof found in \cite{DavSem} relating to the Mumford-Shah problem. 

\begin{lem}\label{DiniVan}
 If $g^j(x)\in L^2(Q_2^+)$, for $j=1,2,...,N$, are $N$ functions satisfying, for every $t\in (0,1)$, 
 $\|g^j(\cdot,t)\|_{L^2(Q_2'(0)\times\{x_n=t\})}\le C_0\sqrt{t}$. Then there is a sequence of cubes $Q_{r_k}(y^k)$, $y^k\in \{x_n=0\}$ such that
 as $k\to \infty$, $r_k\to 0$, and
 $$
 H^j(r_k,y^k)=_{df}\frac{1}{|Q_{r_k}^+(y^k)|}\int_{Q_{r_k}}|g^j(x)|^2dx\to 0,
 $$
 for all $j=1,2,...,N$.
\end{lem}
 \textsl{Proof:} We make the following calculation, where we use the definition of $H^j$ in the first step,
 $$
 \int_{y'\in Q_1'\times\{x_n=0\}}\int_{t=0}^1 \frac{1}{t}H^j(t,y)dtdy'=
 $$
 \begin{equation}\label{BoringStuff}
 = \int_{y'\in Q_1'\times\{x_n=0\}} \left(\int_{t=0}^1 \left(\frac{1}{t^{n+1}}\int_{Q_t^+(y')} |g^j(x)|^2dx\right)dt\right)dy'=
 \end{equation}
 $$
 = \int_{y'\in Q_1'\times\{x_n=0\}} \left(\int_{t=0}^1 \left(\frac{1}{t^{n+1}}\int_{x_n=0}^t\int_{x'\in Q_t'(y')} |g^j(x',x_n)|^2dx\right)dt\right)dy'=
 $$
 $$
 =\int_{x_n=0}^1 \int_{t=x_n}^1 \frac{1}{t^{n+1}}\int_{y'\in Q_1'\times\{x_n=0\}} \int_{x'\in Q_t'(y')} |g^j(x',x_n)|^2dx' dy' dt dx_n,
 $$
 where we used Fubini's Theorem on the two middle integrals in the last equality.
 
 We may continue to estimate the right side of (\ref{BoringStuff}) by noticing that integrating over $x'\in Q'_t(y')$ and then over $y'\in Q_1'(0)$
 may be estimated by $t^{n-1}$ times an integration of $y'$ over $Q'_{1+t}\subset Q_2'$ for small $t$. This gives that (\ref{BoringStuff}) can be estimated 
 from above by
 $$
 \int_{x_n=0}^1 \left(\int_{t=x_n}^1\frac{1}{t^{n+1}}\left(t^{n-1}\int_{x'\in Q_2'(0)} |g^j(x',x_n)|^2dx'\right)dt\right)dx_n=
 $$
 \begin{equation}\label{StillBoring}
 =\int_{x_n=0}^1 \left(\int_{t=x_n}^1\frac{1}{t^{2}}dt\left(\int_{x'\in Q_2'(0)} |g^j(x',x_n)|^2dx'\right)\right)dx_n=
 \end{equation}
 $$
 =\int_{x_n=0}^1 \left( \frac{1}{x_n}-1\right)\left(\int_{x'\in Q_2'(0)} |g^j(x',x_n)|^2dx'\right)dx_n\le C_0^2,
 $$
 where we used the assumption that $\|g^i(\cdot, t)\|_{L^2}^2\le C_0^2t$ in the last inequality.
 
 Putting the estimates (\ref{BoringStuff}) and (\ref{StillBoring}) together we have shown that 
 \begin{equation}\label{AllTheBoringTogether}
  \int_{y'\in Q_1'\times\{x_n=0\}}\int_{t=0}^1 \frac{1}{t}H^j(t,y)dtdy'\le C_0^2.
 \end{equation}

 Observe that this implies that, for every $\epsilon>0$, if we define the the set
 $$
 K_\epsilon=\left\{ (t,y);\; y\in Q_1'\times\{x_n=0\},\; t\in (0,1), \sup_{j}\left(H^j(t,y)\right)\ge \epsilon \right\}.
 $$
 then 
 \begin{equation}\label{Triviality}
 \epsilon\chi_{K_\epsilon}(t,y)\le H^j(t,y).
 \end{equation}
 
 From (\ref{Triviality}) and (\ref{AllTheBoringTogether}) it follows that
 \begin{equation}\label{Handthatfeeds}
 \int_{y\in Q_1'\times\{x_n=0\}}\int_{t=0}^1 \frac{1}{t} \chi_{K_\epsilon}dtdy\le \frac{C_0^2}{\epsilon}. 
 \end{equation}
 This implies that, for every $\delta>0$,
 \begin{equation}\label{Ivarsoga}
  \left\{ (t,y);\; y\in Q_1',\; 0<t<\delta\right\}\setminus K_\epsilon \ne \emptyset.
 \end{equation}
Since if (\ref{Ivarsoga}) was not true then $\chi_{K_\epsilon}(t,y)=1$ for $0<t<\delta$ which would imply that 
$$
\int_{y\in Q_1'\times\{x_n=0\}}\int_{t=0}^1 \frac{1}{t} \chi_{K_\epsilon}dtdy\ge \int_{y\in Q_1'\times\{x_n=0\}}\int_{t=0}^\delta \frac{1}{t} dtdy,
$$
the right integral diverges which would contradict (\ref{Handthatfeeds}); thus (\ref{Ivarsoga}) has to hold.

We may therefore for every $\epsilon_j>0$, say $\epsilon_j=\frac{1}{j}$, and $\delta_k>0$, say $\delta_k=\frac{1}{k}$
find a $(t_{j,k},y_{j,k})\notin K_{\epsilon_j}$ and $0<t_{j,k}<\delta_k$. The diagonal sequence $(r_k,y^k)=(t_{k,k},y_{k,k})$
has the desired property of the Lemma. \qed

\section{Blow-ups of Global Solutions.}

In this section we prove that if $u$ is a global solution in $\R^n$ and $\overline{\Omega_u}$ is a half space, $\{x_n\ge 0\}$, then
$u$ has a blow-up that is one dimensional. The proof is based on Lemma \ref{DiniVan}, in particular we will show that 
the tangential derivatives $\frac{\partial u}{\partial x_i}$, $i=1,2,...,n-1$ satisfies the assumptions in Lemma \ref{DiniVan}.
This implies that there is a blow-up whose $x_j$ derivatives, $j=1,2,...,n-1$, all vanish.

The argument is quite delicate and will partly be done in the space $H^{-\frac{1}{2}}(\partial D)$, for a Lipschitz domain $D$. 
The space $H^{-\frac{1}{2}}(\partial D)$ is the dual space of $H^{\frac{1}{2}}(\partial D)$ - that is the space of traces of 
$W^{1,2}(D)$ functions. In the next lemma we define a trace operator $\gamma$ that gives a $H^{-\frac{1}{2}}(\partial D)$
functional for every vector field $\v$ with divergence in $L^2$. The proof is standard.

\begin{lem}\label{trace}
 Let $D$ be a bounded $C^{0,1}$ domain, $\v\in L^2(D)$ be a vector field with $L^2$ divergence $\d(\v)\in L^2(D)$. Then there is
 a unique $\gamma(\v)\in H^{-\frac{1}{2}}(\partial D)$ such that $\gamma(\v)=\v\big|_{\partial D}$ for each
 $\v\in C^\infty(D)$.
 
 Furthermore:
 $$
 \|\gamma(\v)\|_{H^{-\frac{1}{2}}(\partial D)}\le \|\v\|_{L^2}+\|\d(\v)\|_{L^2}.
 $$
\end{lem}
{\sl Proof:} We define the pairing $\langle \gamma(\v), \phi\rangle$ for $\phi\in H^{\frac{1}{2}}(\partial D)$ according to
\begin{equation}\label{defgamma}
\langle \gamma(\v),\phi \rangle=\int_D \left( \v\cdot \nabla \tilde{\phi}+\d(\v)\tilde{\phi}\right),
\end{equation}
where $\tilde{\phi}$ is the unique $W^{1,2}(D)$ harmonic function satisfying $\tilde{\phi}=\phi$ on $\partial D$.

Then
$$
\left|\langle \gamma(\v),\phi \rangle\right|\le \left( \|\v\|_{L^2(D)}+\|\d(\v)\|_{L^2(D)}\right)\|\phi\|_{H^{\frac{1}{2}}(\partial D)}.
$$
So $\gamma(\v)$ is a bounded functional on $H^{\frac{1}{2}}(\partial D)$ and thus, by Riesz representation Theorem, representable by a 
function in $H^{-1/2}(\partial D)$. 

That $\gamma$ equals the restriction on $C^{\infty}$ vector fields follows from a simple integration by parts in equation (\ref{defgamma}).\qed 

Next we make sure that the tangential derivatives of $u$ is controlled in a way that makes Lemma \ref{DiniVan} applicable. 

\begin{lem}\label{w22estimates}
 Let $u\in W^{1,p}_{loc}(\R^n)$ be a solution to (\ref{main}) in $\R^n$ and $\overline{\Omega}=\R^n_+$. Assume that $u$ satisfies the following 
 growth condition
 \begin{equation}\label{nablagrowth}
 |\nabla u|\le C_0 |x_n|^{\frac{1}{p-1}}
 \end{equation}
 for some $C_0$.
 
 Then, for each $i=1,2,...,n$ and $j=1,2,...,n-1$,
 $$
 \int_{B_1(0)\cap\{0<x_n<t\}}|\nabla u|^{2(p-2)}\left|\frac{\partial^2 u}{\partial x_i \partial x_j}\right|^2\le C_1t,
 $$
 Where $C_1$ depend only on $n$, $p$ and $C_0$.
\end{lem}
{\sl Proof:} Fix a direction $j=1,2,...,n-1$ and let $\partial_h$ be the difference 
quotient in the $j-$direction: $\partial_h f(x)=\frac{f(x+h e_j)-f(x)}{h}$. Then for any function $\phi\in C^\infty_c(B_1(0))$
we may use $\partial_{-h}(\phi^2\partial_h u)$ as a test function in the weak formulation of (\ref{main}) and derive,
$$
\int_{\R^n}\partial_{-h}(\phi^2\partial_h u)=-\int_{\R^n}\nabla \left(\partial_{-h}(\phi^2\partial_h u)\right)\cdot \left( |\nabla u|^{p-2}\nabla u\right)=
$$
$$
=\int_{\R^n}\nabla \left((\phi^2\partial_h u)\right)\cdot \partial_h\left( |\nabla u|^{p-2}\nabla u\right).
$$
The right side is $0$ since it involves integrating a function minus a translate of the function. So a standard rearrangement of terms implies that, 
where we use $I$ for the identity matrix,
$$
\int_{\R^n}\phi^2\nabla \frac{\partial u}{\partial x_j}\cdot \left( (p-2)|\nabla u|^{p-4}\nabla u\oplus \nabla u+|\nabla u|^{p-2}I \right)\nabla \frac{\partial u}{\partial x_j}
\le 
$$
$$
\le 2\int_{\R^n}\phi \frac{\partial u}{\partial x_j} \nabla \phi \cdot \left( (p-2)|\nabla u|^{p-4}\nabla u\oplus \nabla u+|\nabla u|^{p-2}I \right)\nabla \frac{\partial u}{\partial x_j},
$$
which implies that
$$
\int_{\R^n}\phi^2\nabla \frac{\partial u}{\partial x_j}\cdot \left( (p-2)|\nabla u|^{p-4}\nabla u\oplus \nabla u+|\nabla u|^{p-2}I \right)\nabla \frac{\partial u}{\partial x_j}\le
$$
$$
\le C\int_{B_2}|\nabla \phi|^2|\nabla u|^p.
$$
Making a standard choice of $\phi$ and noticing that $(p-2)|\nabla u|^{p-4}\nabla u\oplus \nabla u+|\nabla u|^{p-2}I$ is 
comparable to $|\nabla u|^{p-2}I$ implies that
$$
\int_{B_1(0)\cap\{0<x_n<t\}}|\nabla u|^{p-2}\left|\nabla \frac{\partial u}{\partial x_j}\right|^2\le \frac{C}{t}\int_{B_2(0)\cap\{0<x_n<2t\}}|\nabla u|^p\le
Ct^{\frac{1}{p-1}},
$$
where we used (\ref{nablagrowth}). Using (\ref{nablagrowth}) again gives the result.\qed

\begin{prop}\label{flatnessprop}
 Assume that $u$ is a solution to (\ref{main}) in $\R^n$. Assume furthermore that $\Omega=\R^n_+$ and that
 $$
 \sup_{B_R(0)}|u|\le CR^{\frac{p}{p-1}}.
 $$
 Then there exists an $x^0\in \Gamma=\{x_n=0\}$ and a sequence $r_j\to 0$ such that 
 $$
 \lim_{j\to \infty}\frac{u(r_j x+x^0)}{r_j^{\frac{p}{p-1}}}= \frac{p-1}{p}(x_n^+)^{\frac{p}{p-1}}
 $$
 strongly in $W^{1,p}_{loc}(\R^n)$.
\end{prop}
{\sl Proof:} Taking a derivative, which is well defined in the weak sense, of
$$
\d\left( |\nabla u|^{p-2}\nabla u\right)=\chi_{\{u>0\}}
$$
with respect to $x_i$, $i=1,2,,...,n$, shows that
\begin{equation}\label{derivativeofplap}
\d\left( \left[(p-2)|\nabla u|^{p-4}\nabla u\otimes \nabla u +|\nabla u|^{p-2}I\right]\nabla \v^i\right)=0 \quad \textrm{ in }\R^n_+
\end{equation}
where $\v^i=\frac{\partial u}{\partial x_i}$ and $\nabla u \otimes \nabla u$ is the matrix with entities $u_iu_j$.

With the notation
$$
\w^i=\left[(p-2)|\nabla u|^{p-2}\nabla u\otimes \nabla u +|\nabla u|^{p-2}I\right]\nabla \v^i
$$
it follows, from Lemma \ref{w22estimates}, that $\w^i\in L^2_{\textrm{loc}}(\R^n_+)$ and, from (\ref{derivativeofplap}), that $\w^i$ is divergence 
free in $\R^n_+$.

We can thus conclude from Lemma \ref{trace} that $\w^i$ has a trace 
$$
\gamma_t(\w^i)\in H^{-\frac{1}{2}}_{loc}(\partial \left(\R^n_+\cap \{x_n>t\}\right)).
$$

First we claim that $\gamma_0(\w^i)=0$ for $i=1,2,...,n-1$. To see this we notice that for any test-function $\phi\in C^\infty_c(B_1(x^0))$
we have
$$
\int_{B_1(x^0)}\nabla \phi \cdot \w^i+\phi \d(\w^i)dx=
\int_{B_1(x^0)}\nabla \phi\cdot \left(\frac{\partial}{\partial x_i}  |\nabla u|^{p-2}\nabla u\right)dx=
$$
$$
=-\int_{B_1(x^0)}\nabla \frac{\partial \phi}{\partial x_i}\cdot \left( |\nabla u|^{p-2}\nabla u\right)dx=
$$
$$
=\int_{B_1(x^0)\cap \{x_n>0\}} \frac{\partial \phi}{\partial x_i}dx+\int_{\{x_n=0\}\cap B_1}\phi |\nabla u|^{p-2}\frac{\partial u}{\partial x_n}dx=0,
$$
where we use an integration by parts in several of the steps and the definition of $\w^i$, and that $\d(\w^i)=0$, in the first equality
and that $|\nabla u|^{p-2}\frac{\partial u}{\partial x_n}=0$ on $\{x_n=0\}$, since solutions to the obstacle problem are $C^{1,\alpha}$, in the last equality.

Similarly, for $i=n$, we get
$$
\int_{B_1(x^0)}\nabla \phi \cdot \w^n+\phi \d(\w^n)dx=
\int_{B_1(x^0)}\nabla \phi\cdot \left(\frac{\partial}{\partial x_n}  |\nabla u|^{p-2}\nabla u\right)dx=
$$
$$
=-\int_{B_1(x^0)}\nabla \frac{\partial \phi}{\partial x_n}\cdot \left( |\nabla u|^{p-2}\nabla u\right)dx=
$$
$$
=\int_{B_1(x^0)\cap \{x_n>0\}} \frac{\partial \phi}{\partial x_i}dx=-\int_{B_1(x^0)\cap \{x_n=0\}}\phi(x)dx,
$$
which implies that $\gamma_0(\w^n)=1$. An approximation argument together with the fact that the embedding $C^\infty\subset W^{1,2}$
is dense proves that the above deduction holds for $\phi\in W^{1,2}_0(B_1(x^0))$.

Using that, for $i=1,2,..., n-1$, 
$$
\int_{B_1(x^0)}\nabla \phi \cdot \w^i dx=\int_{B^+_1(x^0)}\nabla \phi \cdot \w^i dx=
$$
$$
=\int_{\partial B_1^+}\phi(x)\gamma_0(\w^i)dx-\int_{B_1^+(0)}\phi(x)\textrm{div}(\w^i(x))dx=0
$$
for any $\phi\in W^{1,2}_0(B_1(x^0))$ we can conclude that
\begin{equation}\label{MeAllaJobb}
\left|\int_{B_1(x^0)\cap\{x_n>t\}}\nabla \phi\cdot \w^i\right|=\left|\int_{B_1(x^0)\cap\{0<x_n<t\}}\nabla \phi\cdot \w^i\right|\le
\end{equation}
$$
\le \|\w^i\|_{L^2(B_1(x^0)\cap \{0<x_n<t\})}\|\nabla \phi\|_{L^2(B_1(x^0)\cap \{0<x_n<t\})}.
$$

From Lemma \ref{w22estimates} and (\ref{MeAllaJobb}) we can conclude that 
\begin{equation}\label{WithLStar}
\langle \gamma_t(\w^i), \phi \rangle\le Ct^{\frac{1}{2}}\|\nabla \phi\|_{L^2(B_1(x^0)\cap \{0<x_n<t\})}.
\end{equation}

Similarly, for $i=n$ we can conclude that
\begin{equation}\label{WithLTwoStar}
\langle \gamma_t(\w^n)-1, \phi \rangle\le Ct^{\frac{1}{2}}\|\nabla \phi\|_{L^2(B_1(x^0)\cap \{0<x_n<t\})}.
\end{equation}

Let us remind ourselves, before we continue with the proof, that we may characterize functions in 
$H^{-\frac{1}{2}}(\R^{n-1})$ by $L^2-$functions. In particular, for any $\kappa\in H^{-\frac{1}{2}}$
there is a $g\in L^2$ such that if $\phi\in H^{\frac{1}{2}}$ then
$$
\langle \kappa,\phi \rangle=\int g(x) \Lambda^{\frac{1}{2}} \phi dx
$$
where $\Lambda^{\frac{1}{2}}$ is defined
$$
\Lambda^{\frac{1}{2}} \phi= \mathcal{F}^{-1}\left( (1+|\xi|^2)^{\frac{1}{4}}\mathcal{F}(\phi)\right)
$$
where $\mathcal{F}$ denotes the Fourier transform.

We may thus identify $\gamma_t(\w^i)$ and $\gamma_t(\w^n)-1$ with $L^2-$functions $g^j(x',t)$ satisfying,
in view of the inequalities (\ref{WithLStar}) and (\ref{WithLTwoStar}),
$$
\int_{Q_2'}|g^j(x',t)|^2dx\le C_0 t
$$
for some constant $C_0$. Using Lemma \ref{DiniVan} we can find a sub-sequence $y^k\in Q_1(0)\cap\{x_n=0\}$ and
$r_k\to 0$ such that $\frac{1}{r_k^{\frac{n+1}{2}}}\|g^j\|_{L^2(Q_{r_k}^+(y^k))}\to 0$. 

We claim that by, possibly considering another sequence, $r_k\to 0$ 
we may assume that 
\begin{equation}\label{LocalHandling}
\frac{1}{r_k^{\frac{n+1}{2}}}\|g^j\|_{L^2(Q_{R r_k}^+(y^k))}\to 0
\end{equation}
for any $R>0$. This is easy to see: if
$$
\frac{1}{r_k^{\frac{n+1}{2}}}\|g^j\|_{L^2(Q_{r_k}^+(y^k))}=\epsilon_k\to 0
$$
then, with $s_k=\epsilon_k^{\frac{1}{n+1}}r_k$, a simple calculation gives
$$
\frac{1}{r_k^{\frac{n+1}{2}}}\|g^j\|_{L^2(Q_{R s_k}^+(y^k))}<\epsilon_k^{1/2}\to 0,
$$
for every $R<\frac{1}{\epsilon_k^{1/(n+1)}}\to \infty$. We may assume that $r_k$ already satisfies (\ref{LocalHandling}).

Since $\sup_{B_R}u^k\le C R^{p/(p-1)}$ by Lemma \ref{RegularityPharmObst} we may conclude that, for some subsequence,
$$
u^k(x)=\frac{u(r_k x+ y^k)}{r_k^{\frac{p}{p-1}}}\to u^0(x).
$$

We will show that the the limit $u^0$ satisfies  
\begin{equation}\label{U0Eq1}
 \Delta_p u^0=\chi_{\{u^0>0\}} \quad \textrm{ in }\R^n 
\end{equation}
\begin{equation}\label{U0Eq2}
 \sup_{B_R(0)}|u^0|\le C\left( 1+R^{\frac{p}{p-1}}\right) \quad \textrm{ for every }R>0 
\end{equation}
\begin{equation}\label{U0Eq3}
 \gamma_t(\w^j_{u^0})=0 \quad \textrm{ for }j=1,2,...,n-1 \textrm{ and all }t>0,
\end{equation}
\begin{equation}\label{U0Eq4}
 \gamma_t(\w^n_{u^0})=1 \textrm{ for all }t>0,  
\end{equation}
where $\w^j_{u^0}=\left[(p-2)|\nabla u^0|^{p-2}\nabla u^0\otimes \nabla u^0 +|\nabla u^0|^{p-2}I\right]\nabla \frac{\partial u^0}{\partial x_j}$.
That (\ref{U0Eq1}) holds is trivial, (\ref{U0Eq2}) follows from Lemma \ref{RegularityPharmObst} and (\ref{U0Eq3}) and (\ref{U0Eq4}) follows from 
the fact that $g^j\to 0$ by (\ref{LocalHandling}).

But since $\w^j_{u^0}\in L^2$ we can conclude, from (\ref{U0Eq3}), that $\w^j_{u^0}=0$ almost everywhere, for $j=1,2,...,n-1$. 
Similarly, by (\ref{U0Eq4}), $\w^n_{u^0}=e_n$ almost everywhere.

Since $\w^j_{u^0}=\frac{\partial }{\partial x_j}|\nabla u^0|^{p-2}\nabla u^0$, for $j=1,2,...,n$, and $\w^j_{u^0}=0$ in $x_n<0$ we can conclude that
$|\nabla u^0|^{p-2}\frac{\partial u^0}{\partial x_n}=x_n$.

It directly follows, from the definition of $\w^j_{u^0}$, that 
\begin{equation}\label{HubbaBubba}
\nabla'\cdot \left( |\nabla u^0|^{p-2}\nabla' u^0\right)=0 
\end{equation}
for each $x_n>0$ and that $|\nabla u^0|\ge x_n^{\frac{1}{p-1}}$
for $x_n>0$. 

Since $u^0(x',0)=0$ it follows from Lemma \ref{RegularityPharmObst} that $u^0(x', x_n)$ is uniformly bounded in $x'$ for each fixed $x_n$. From 
(\ref{HubbaBubba}) it also follows that $\frac{\partial u^0(x', x_n)}{\partial x_j}$ is constant in $x'$ for each fixed $x_n$;
this together with the boundedness of $u^0(x',x_n)$ implies that $u^0(x',x_n)$ is constant in $x'$.
That is $u^0(x)=u^0(x_n)$, and 
\begin{equation}\label{SomeODE}
\Delta_p u^0(x)=\frac{\partial }{\partial x_n}\left(\left|\frac{\partial u^0(x_n)}{\partial x_n}\right|^{p-2}\frac{\partial u^0(x_n)}{\partial x_n} \right)=0
\end{equation}
Solving (\ref{SomeODE}) gives that $u^0(x_n)=\frac{p-1}{p}x_n^{\frac{p}{p-1}}$ in the set $\{x_n>0\}$. \qed

\section{Linearization.}\label{sec:lin}

In this section we show that we may linearize solutions to the p-harmonic obstacle problem around points 
where the solution is close to the ground-state solution $\frac{p-1}{p}(x_n)_+^{\frac{p}{p-1}}$. Before we do that 
we fix some notation that we will use in this section. Notice that if $u$ is a solution to the p-harmonic obstacle 
problem then $v_i=\frac{\partial u}{\partial x_i}$ will be a weak solution to the following linear equation
$$
\d(A(x)\nabla v_i)=0 \quad\quad\quad\textrm{ in } \Omega=\{u>0\},
$$
where $A(x)=(p-2)|\nabla u|^{p-4}\nabla u \otimes \nabla u+|\nabla u|^{p-2} I$.

We will use $L^2_{A}(\Omega)$ for the Hilbert space with the norm 
$$
\|v\|_{L^2_{A}(\Omega)}=\left( \int_{\Omega}|A(x)||v|^2 dx\right)^{1/2}.
$$
At times we will be somewhat informal and say that $\nabla v\in L^2_{A}(\Omega)$ if
$$
\left( \int_{\Omega}\langle \nabla v \cdot A(x), \nabla v\rangle dx\right)^{1/2}<\infty,
$$
which is justified since $|\langle \nabla v \cdot A(x), \nabla v\rangle|$ is comparable to $|A||\nabla v|^2$ 
up to a multiplicative constant that only depends on $(p-2)$. Thus we can form the Hilbert space $W^{1,2}_{A}$ with norm
$\|v\|_{L^2_{A}}+\|\nabla v\|_{L^2_{A}}$.

In this section we will assume that we have a sequence of solutions $u^j$ in $B_2(0)$ in a normalized coordinate system defined as follows.
\begin{definition}
Let $u\in W^{1,p}(B_1(0))$ be a solution to (\ref{main}). Then we say that the coordinate system is normalized (with respect to
$u$) if
$$
\left\| \nabla' \left(u-\frac{p-1}{p}x_n^{\frac{p}{p-1}}\right) \right\|_{W^{1,2}_{A}(\Omega\cap B_1)}\le 
$$
$$
\le \left\|  \nabla' \left(u-\frac{p-1}{p}(\nu\cdot x + \gamma)^{\frac{p}{p-1}} \right)\right\|_{W^{1,2}_{A}(\Omega\cap B_1)}
$$
and $\gamma$ is chosen so that 
$$
\left\| \left(u-\frac{p-1}{p}x_n^{\frac{p}{p-1}}\right) \right\|_{L^2_{A}(\Omega\cap B_1)}\le 
$$
$$
\le \left\| \left(u-\frac{p-1}{p}(x_n + \gamma)_+^{\frac{p}{p-1}} \right)\right\|_{L^2_{A}(\Omega\cap B_1)}
$$
for any unit vector $\nu$ and constant $\gamma$, where 
$$
A=(p-2)|\nabla u|^{p-4}\nabla u \otimes \nabla u+|\nabla u|^{p-2} I,
$$
and $\nabla'_\nu=\nabla-\nu(\u\cdot \nabla)$ is the gradient on the subspace orthogonal to $\nu$.
\end{definition}
Notice that if $u$ is a solution to (\ref{main}) in some set we can easily chose a normalized coordinate system by making a translation and rotation of the 
coordinate system. Also, the term $\frac{p-1}{p}x_n^{\frac{p}{p-1}}$ is redundant in the first set of inequalities since it disappears under the 
operation of $\nabla'$, we have included the term as an indication of the heuristic idea that the normalized coordinate system 
is the coordinate system where $u$ is closest to a solution only depending on the $x_n$ coordinate.

For the rest of this section we will assume that we assume that $u^j$ is a sequence of solutions to (\ref{main}) such that
$$
\left\| \nabla'\left(u^j-\frac{p-1}{p}(x_n)^{\frac{p}{p-1}}\right) \right\|_{W^{1,2}_{A_j}(\Omega_j\cap B_1)}\to 0,
$$
we may think of $u^j(x)=\frac{u(r_jx)}{r_j^{p/(p-1)}}$ for some function p-harmonic function $u(x)$ whose blow-up at the origin is 
$\frac{p-1}{p}(x_n)_+^{\frac{p}{p-1}}$. 

In the rest of this section we will also denote by  $v^j_i$, for $i=1,2,...,n-1$, the functions
\begin{equation}\label{DefinitionofVij}
v^j_i(x)=\frac{1}{\delta_{i,j}}\frac{\partial u^j}{\partial x_i},
\end{equation}
where $\delta_{i,j}>0$ is chosen such that $\|v^j_i\|_{L^2(\Omega_j)}=1$. This is only well defined if 
$\left\|\frac{\partial u^j}{\partial x_i}\right\|_{L^2(\Omega_j)}\ne 0$, but we may throw out the terms in the sequence where this is not satisfied. 
In the very special case when $\left\|\frac{\partial u^j}{\partial x_i}\right\|_{L^2(\Omega_j)}=0$ for all but finitely many $j$ we can throw out the 
terms where $\left\|\frac{\partial u^j}{\partial x_i}\right\|_{L^2(\Omega_j)}\ne 0$ and consider the sequence $u^j$ as a sequence of solutions in 
$\R^{n-1}$. Therefore there is no loss of generality to assume that all $v^j_i$ are well defined.

\begin{prop}\label{StrongTangConv}
 Assume that $p>2$. The sequence $v^j_i$ defined by (\ref{DefinitionofVij}), for $i=1,2,...,n-1$, converges to $v^0_i$ where $v^0_i$ satisfies
 $$
 \d\left(|x_n|^{\frac{p-2}{p-1}} \nabla v^0_i\right) = 0 \quad \textrm{ in } \{ x_n>0\}\cap B_1(0)
 $$
 and 
 $v^0_i=0$ on $\{x_n=0\}\cap B_1(0)$. The convergence $v^j_i\to v^0_i$ is in the sense that for any $\epsilon>0$
 \begin{equation}\label{strongtangential}
 v^j_i\to v^0_i \quad\quad \textrm{ in } C^{1,\alpha}(\{x_n>\epsilon\}) 
 \end{equation}
 and $\limsup_{j\to \infty} \|v^j_i\|_{L^2(B_{1/2})}=\|v^0_i\|_{L^2(B_{1/2}(0))}$.
\end{prop}
\textsl{Proof:} That $\nabla u^j \to \nabla u^0$ in $C^\beta$ and weakly in $L^p$, for a sub-sequence, follows from the $C^{1,\alpha}$ regularity
of $u^j$ together with the Arzela-Ascoli theorem and weak compactness in $L^p$ spaces. Moreover, by Lemma \ref{L2} and and $C^\beta$ convergence of $A_j$,
it follows that $A_j\cdot \nabla v^j_i\to A_0\cdot \nabla v^0_i$ weakly and that 
$$
0=\int_{B_1}\langle \nabla \phi \cdot A_j, \nabla v^j_i\rangle=\int_{B_1}\langle \nabla \phi \cdot A_0, \nabla v^0_i\rangle.
$$

The only thing that remains to prove is the strong convergence of $v^j_i$ specified at the end of the Lemma. We 
will prove (\ref{strongtangential}) by means of a hole filling argument.

In order to prove the strong convergence in $L^2$ we fix any $\epsilon>0$. Then, since $u^j\to \frac{p-1}{p}(x_n)_+^{\frac{p}{p-1}}$ in $C^{1,\alpha}$
and $u^j$ is non-degenerate by Lemma \ref{RegularityPharmObst} (in particular (\ref{OptGrowth})), 
if $j$ is large enough it follows that $\partial \{u^j>0\}\cap B_1 \subset \{|x_n|<\epsilon\}$. We can thus, for any $\delta>0$, 
assume, if $j$ is large enough, that the entire free boundary is contained in the strip $|x_n|<\delta$.

Let us define the following cubes
$$
Q_k=\{x\in B_1(0);\; x'\in B'_{1/2+(2^{k+1}-1)\delta}(0) \textrm{ and } |x_n|<(2^{k+1}-1)\delta\}
$$
then $Q_k\in B_{3/4}$ for $k\le \frac{\ln(\delta)}{\ln(2)}$. Furthermore we have 
\begin{enumerate}
 \item There exists cut-off functions $\Psi_k\in C^\infty$ such that $\psi^k=1$ in $Q_k$, $\psi^k=0$ in $Q_{k+1}^c$ and $|\nabla \psi_k|\le C \frac{1}{2^k \delta}$.
 \item $|u^j|\le C 2^{\frac{pk}{p-1}}\delta^{\frac{p}{p-1}}$ and $|\nabla u^j|\le C 2^{\frac{k}{p-1}}\delta^{\frac{1}{p-1}}$ for 
 some constant $C$ (depending only on $n$ and $p$) in $Q_k$. This follows from Lemma \ref{RegularityPharmObst} and that $\Gamma_u\in \{|x_n|<\delta\}$.
\end{enumerate}

\begin{lem}\label{holefilling}
 There exists a constant $\tau<1$, depending only on $n$ and $p\ge 2$, such that for every $k\ge 0$
 $$
 \int_{Q_k}|v^j_i|^2dx \le \tau \int_{Q_{k+1}}|v^j_i|^2dx.
 $$
\end{lem}
\textsl{Proof of Lemma \ref{holefilling}:} We use $\phi=\psi_k^2 (v^j_i)^2$ as a test function in the variational integral for 
$u^j$, this is well defined since $v^j_i\in C^{\alpha}\cap W^{1,2}_{A_j}$ by Lemma \ref{RegularityPharmObst} and Lemma \ref{L2}. That is
$$
0=\int_{B_1}\left(\nabla \left(\psi_k^2(v_i^j)^2\right)\cdot \nabla u^j |\nabla u^j|^{p-2}+ \left(\psi_k^2(v_i^j)^2\right)\right)dx=
$$
$$
=\int_{B_1}\Big(\left( 2\psi_k(v^j_i)^2\nabla \psi_k\cdot \nabla u^j |\nabla u^j|^{p-2}\right)+
$$
$$
+\left(2\psi_k^2 v^j_i\nabla v^j_i\cdot \nabla u^j |\nabla u^j|^{p-2} \right)+ \left(\psi_k^2(v_i^j)^2\right)\Big)dx.
$$
Rearranging terms we can deduce that 
$$
\int_{B_1}\psi_k^2(v_i^j)^2dx\le \left|\int_{B_1}\left( 2\psi_k(v^j_i)^2\nabla \psi_k\cdot \nabla u^j |\nabla u^j|^{p-2}\right)dx \right|+
$$
\begin{equation}\label{Some1Star}
+\left|\int_{B_1}\left(2\psi_k^2 v^j_i\nabla v^j_i\cdot \nabla u^j |\nabla u^j|^{p-2} \right)dx \right|=I_1+I_2.
\end{equation}
We will estimate $I_1$ and $I_2$ independently.

In order to estimate $I_1$ we use that $\nabla \psi_k$ is supported in $Q_{k+1}\setminus Q_k$ where 
$|\nabla u^j|^{p-1}\le C  2^k\delta\le \frac{C}{\|\nabla \psi_k \|}$.
Therefore 
\begin{equation}\label{Some2Star}
I_1=\left|\int_{Q_{k+1}\setminus Q_k}\left( 2\psi_k(v^j_i)^2\nabla \psi_k\cdot \nabla u^j |\nabla u^j|^{p-2}\right)dx \right|\le 
C \int_{Q_{k+1}\setminus Q_k}\psi_k(v^j_i)^2dx.
\end{equation}

In order to estimate $I_2$ we use a small trick and notice that 
$\nabla u^j |\nabla u^j|^{p-2}=\frac{1}{p-1}A_j\cdot \nabla u^j$. We may therefore estimate $I_2$
$$
I_2=\frac{1}{p-1}\left|\int_{B_1}\left(2\psi_k^2 v^j_i\langle \nabla v^j_i\cdot A_j,\nabla u^j\rangle \right)dx \right|=
$$
\begin{equation}\label{LikeNumbers}
=\frac{1}{p-1}\Bigg|\int_{Q_{k+1}}\left(2\psi_k^2 u^j \langle \nabla v^j_i\cdot A_j,\nabla v^j_i\rangle \right)dx+ 
\end{equation}
$$
+\int_{Q_{k+1}\setminus Q_k}\left(4\psi_k v^j_i u^j \langle \nabla v^j_i\cdot A_j,\nabla \psi_k\rangle \right)dx\Bigg|
$$
where we used an integration by parts and $\textrm{div}(A_j \nabla v^j_i)=0$.

Continuing (\ref{LikeNumbers}) by using the triangle inequality and then H\"older's inequality  it follows that
$$
I_2 \le \frac{1}{p-1}\left|\int_{B_1}\left(2\psi_k^2 u^j \langle \nabla v^j_i\cdot A_j,\nabla v^j_i\rangle \right)dx\right|+ 
$$
$$
+\frac{1}{p-1}\left|\int_{Q_{k+1}\setminus Q_k}\left(4\psi_k v^j_i u^j \langle \nabla v^j_i\cdot A_j,\nabla \psi_k\rangle \right)dx\right|\le
$$
$$
\le \frac{C}{p-1}\Bigg(2^{\frac{pk}{p-1}}\delta^{\frac{p}{p-1}}\left|\int_{Q_{k+1}}\psi_k^2\langle \nabla v^j_i\cdot A_j,\nabla v^j_i\rangle dx \right|+
$$
$$
+\left(2^{\frac{pk}{p-1}}\delta^{\frac{p}{p-1}}\int_{Q_{k+1}\setminus Q_k} \psi_k^2\langle \nabla v^j_i\cdot A_j,\nabla v^j_i\rangle dx\right)^{\frac{1}{2}}  \left(\int_{Q_{k+1}\setminus Q_k}(v^j_i)^2dx\right)^{\frac{1}{2}}
\Bigg),
$$
where we also used that $|u^j|\le C 2^{\frac{pk}{p-1}}\delta^{\frac{p}{p-1}}$, $|\nabla \psi_k|\le \frac{C}{2^k \delta}$ and 
$|\nabla u|\le C2^{\frac{k}{p-1}}\delta^{\frac{1}{p-1}}$, see point 1 and 2 in the statement just before this Lemma, in $Q_{k+1}$ in the last step.

From Lemma \ref{L2} it follows that, since $p\ge 2$,
$$
\int_{Q_{k+1}}\psi_k^2\langle \nabla v^j_i\cdot A_j,\nabla v^j_i\rangle dx\le 
\frac{C\sup_{Q_{k+1}}|\nabla u|^{p-2}}{2^{2k}\delta^2} \int_{Q_{k+1}\setminus Q_k} (v^j_i)^2dx\le
$$
$$
\le C 2^{-\frac{pk}{p-1}}\delta^{-\frac{1}{p-1}}\frac{C\sup_{Q_{k+1}}|\nabla u|^{p-2}}{2^{2k}\delta^2} \int_{Q_{k+1}\setminus Q_k} (v^j_i)^2dx.
$$
We may therefore estimate $I_2$ is according to
\begin{equation}\label{Some3Star}
I_2\le C \int_{Q_{k+1}\setminus Q_k}(v^j_i)^2dx.
\end{equation}

From (\ref{Some1Star}), (\ref{Some2Star}) and (\ref{Some3Star}) we can conclude that 
$$
\int_{Q_k}(v_i^j)^2dx\le C \int_{Q_{k+1}\setminus Q_k} (v_i^j)^2 dx,
$$
the result follows by adding $C\int_{Q_k}(v_i^j)^2dx$ to both sides and dividing by $C+1$. \qed

\textsl{The end of the proof of Proposition \ref{StrongTangConv}:} It is enough to show that for every $\epsilon>0$ there exists a $J_\epsilon$ such that
$j>J_\epsilon$ implies that
$$
\|v^j_i\|_{L^2(B_{1/2})}-\epsilon \le \|v^0_i\|_{L^2(B_{1/2}(0))} \le \|v^j_i\|_{L^2(B_{1/2})}+\epsilon.
$$
Since $v^j_i\to v^0_i$ in $C^{\infty}(B_{1/2}(0)\setminus Q_0)$ for  $\delta>0$ it is enough to show that
there exists an $J_{\epsilon}$ such that $\|v^j_i\|_{L^2(B_{1/2}\cap Q_0)}^2< \epsilon$ for every $j>J_{\epsilon}$. But
by Lemma \ref{holefilling} we know that
\begin{equation}\label{ressurection}
\|v^j_i\|_{L^2(B_{1/2}\cap Q_0)}^2\le \tau^k \|v^j_i\|_{L^2(Q_k)}^2 \le \tau^k
\end{equation}
for every $k\le \frac{\ln(\delta)}{\ln(2)}$. 

We may thus, for every $\epsilon>0$ choose $k_0$ large enough so that $\tau^{k_0}<\epsilon$ and then choose $\delta>0$ so that 
$k_0\le \frac{\ln(\delta)}{\ln(2)}$. Then, if $j$, is large enough we may apply (\ref{ressurection}) and conclude that
$$
\|v^j_i\|_{L^2(B_{1/2}\cap Q_0)}^2\le \tau^{k_0} \|v^j_i\|_{L^2(Q_k)}^2 \le \tau^{k_0}<\epsilon.
$$
This finishes the proof.
\qed

We also need to control the $x_n-$derivative. To that end we prove the following simple lemma.

\begin{lem}\label{xncontrol}
 Given $M$ there exists a modulus of continuity $\sigma$ depending only on $2<p<\infty$ and the dimension $n$ such that 
 \begin{equation}\label{xncontrolEQ}
 \left\|\frac{\partial (u-(p-1)/p (x_n)_+^{p/(p-1)})}{\partial x_n}\right\|_{L^2(\Omega_u\cap B_1)}\le \sigma\left( \|\nabla' u\|_{L^2(B_1)}\right)
 \end{equation}
 for any solution $u$ to (\ref{main}) for which $0\in \Gamma_u$. 
\end{lem}
\textsl{Proof:} We argue by contradiction and assume that $u^j$ is a sequence of solutions satisfying the assumptions
in the Lemma and 
$$
\|\nabla' u^j\|_{L^2(B_1)}\to 0
$$
and 
\begin{equation}\label{Nonsense}
\left\|\frac{\partial (u^j-((p-1)/p) (x_n)_+^{p/(p-1)})}{\partial x_n}\right\|_{L^2(\Omega_u\cap B_1)}\ge \delta
\end{equation}
for some $\delta>0$. 

Since $0\in \Gamma_{u^j}$ it follows from Lemma \ref{RegularityPharmObst} that $\|u^j\|_{L^\infty}(B_{5/4}(0))$ is
bounded and therefore, also by Lemma \ref{RegularityPharmObst}, that $\|u^j\|_{C^{1,\alpha}(B_{6/5}(0))}$ is bounded.

We may therefore choose a subsequence, still denoted by $u^j$, such that $u^j\to u^0$ in $C^{1,\beta}$ where $u^0$ is a solution satisfying 
$\|\nabla' u^0\|_{L^2(B_1)}=0$. In particular, $u^0$ is a one dimensional solution, with $0\in \Gamma_{u^0}$, and thus $u^0=\frac{p-1}{p}x_n^{\frac{p}{p-1}}$ in its support. This 
together with $C^{1,\beta}$ convergence contradicts (\ref{Nonsense}).\qed

\section{Almost Everywhere Uniqueness of Blow-ups.}

We are now ready to prove geometric decay for the tangential derivatives which will lead to $C^{1,\alpha}-$regularity of the free boundary.

\begin{prop}\label{improvement}
 For every $M>0$ there exists constants $\epsilon_0>0$, $\tau<1$ and $s>0$, depending only on $2\le p<\infty$ and $n$, such that
 if $u$ is a solution to (\ref{main}) in a normalized coordinate system such that $\Gamma_u\cap B_{1/4}(0)\ne \emptyset$ 
 and $\|\nabla' u\|_{L^2(B_1(0))}\le\epsilon <\epsilon_0$ then
 there exists a coordinate system and an $0<s<1$ such that
 \begin{equation}\label{SomeCrap}
 \left\| \nabla' \frac{u(sx)}{s^{\frac{p}{p-1}}} \right\|_{L^2(B_1(0))}\le \tau \epsilon.
 \end{equation}
\end{prop}
\textsl{Proof:} We will argue indirectly and assume that $u^j$ is a sequence of solutions as in the Proposition with 
$$
\|\nabla' u^j\|_{L^2(B_1(0))}\le\frac{1}{j}.
$$
Then, by Lemma \ref{xncontrol}, 
$$
\left\|\frac{\partial \left(u-\frac{p-1}{p} x_n^{p/(p-1)}\right)}{\partial x_n}\right\|_{L^2(\Omega_u\cap B_1)}\le \sigma\left( \frac{1}{j}\right)\to 0.
$$
From Proposition \ref{StrongTangConv} we may therefore conclude that, for $i=1,2,...,n-1$,
\begin{equation}\label{New15}
v^j_i=\frac{1}{\delta_{i,j}}\frac{\partial u^j}{\partial x_i}\to v^0_i,
\end{equation}
in $L^2(B_r(0))$ for every $r<1$, where 
$$
\delta_{i,j}=\left\|\frac{\partial u^j}{\partial x_i}\right\|_{L^2(B_1)}
$$
and $v^0_i$ solves
$$
 \d\left(|x_n|^{\frac{p-2}{p-1}} \nabla v^0_i\right) = 0 \quad \textrm{ in } \{ x_n>0\}\cap B_1(0).
$$
By Lemma \ref{eigfunclem} we may write 
\begin{equation}\label{DefOfVi0}
v^0_i(x)=\sum_{k=1}^\infty a_k^i q_j(x),
\end{equation}
where $q_j$ are $\lambda_j-$homogeneous eigenfunctions on the sphere. Also $q_1(x)=x_n^{\frac{1}{p-1}}$.

We will show that $a_1^i=0$ for $i=1,2,...,n-1$. As we will show at the end of the proof; this implies that $v^j_i$ decays faster than $s^{p/(p-1)}$
which will imply (\ref{SomeCrap}) for $u^j$ when $j$ is large. In order to show that $a_1^i=0$ for $i=1,2,...,n-1$ we define 
$$
G(\eta)=\int_{B_1(0)}A_j(x)\left|\nabla'_{\eta} \left(u(x)-\frac{p-1}{p}(\eta\cdot x)_+^{p/(p-1)}\right)\right|^2dx,
$$
on the set of unit vectors $\eta\in \R^n$. Since we assume that the coordinate system is normalized with respect to $u$ it follows that
\begin{equation}\label{EnIsMin}
G(e_n)\le G(\eta),
\end{equation}
for all unit vectors $\eta$. If we let $\eta'=(\eta_1,...,\eta_{n-1},0)$ be a unit vector then $\eta(\theta)=\cos(\theta)e_n+\sin(\theta)\eta'$
will also be a unit vector and (\ref{EnIsMin}) states that $G(\eta(0))\le G(\eta(\theta))$. Taking a derivative with respect to $\theta$ at $\theta=0$
gives 
$$
0=-2\int_{B_1(0)}A_j(x)\frac{\partial u^j}{\partial x_n} \nabla' u^j\cdot \eta' dx=
$$
\begin{equation}\label{ThisISTheShit}
=-2\int_{B_1(0)}A_j(x)\frac{\partial u^j}{\partial x_n} \left(\sum_{i=1}^{n-1}\delta_{i,j}v^j_i(x)\eta_i \right)dx,
\end{equation}
for any unit vector $\eta'$.

Since $\frac{\partial u^j}{\partial x_n}\to (x_n)_+^{1/(p-1)}=q_1(x)$ strongly and $v^j_i\rightharpoonup v_i^0$, defined as in (\ref{DefOfVi0})
we may conclude, by choosing $\eta'=e_i$ dividing by $\delta_{i,j}$ and use that $(x_n)_+^{p/(p-1)}=q_1(x)$, that 
$$
0=\int_{B_1(0)}|x_n|^{\frac{p-2}{p-1}} \sum_{k=1}^\infty a_k^i q_j(x)q_1(x) dx=a_1^i\int_{B_1(0)}|x_n|^{\frac{p-2}{p-1}} |q_1(x)|^2 dx,
$$
where we used that $q_k$ and $q_l$ are orthogonal, by Lemma \ref{eigfunclem}, in the last equality. It follows that $a_1^i=0$ for $i=1,2,...,n-1$.

We have therefore shown that, for any $i=1,2,...,n-1$,
$$
v^j_i\to v^0_i=\sum_{k=2}^\infty a_k^i q_k(x)
$$
strongly (by Proposition \ref{StrongTangConv}) in $L^2(B_{1/2})$ where
$$
\left\| \sum_{k=2}^\infty a_k^i q_k(x) \right\|_{L^2(B_1(0))}\le 1.
$$
Since $q_k$ is homogeneous of order $\lambda_k>\frac{1}{p-1}$ for $k\ge 2$ by Lemma \ref{eigfunclem} it follows that
\begin{equation}\label{skrapigt}
\left\|\frac{v^0_i(sx)}{s^{1/(p-1)}}\right\|_{L^2(B_1)}\le s^{\lambda_2-\frac{1}{p-1}}\|v^0_i\|_{L^2(B_1)}\le s^{\lambda_2-\frac{1}{p-1}}.
\end{equation}

Finally we notice that
$$
\epsilon=\left\| \nabla' u^j(x)\right\|_{L^2(B_1(0))}^2=\sum_{i=1}^{n-1}\delta_{i,j}^2\left\|v^j_i(x)\right\|_{L^2(B_1)}^{2}=\sum_{i=1}^{n-1}\delta_{i,j}^2
$$
so if $\epsilon$ is small enough then
$$
\left\| \nabla' \frac{u^j(sx)}{s^{p/(p-1)}}\right\|_{L^2(B_1(0))}^2=\sum_{i=1}^{n-1}\delta_{i,j}^2\left\|\frac{v^j_i(sx)}{s^{1/(p-1)}}\right\|_{L^2(B_1)}^{2}=
$$
$$
=\sum_{i=1}^{n-1}\delta_{i,j}^2\left\|\frac{v^0_i(sx)}{s^{1/(p-1)}}\right\|_{L^2(B_1)}^{2}+o(\epsilon)\le 
$$
$$
\le \sum_{i=1}^{n-1}s^{\lambda_2-\frac{1}{p-1}}\delta_{i,j}^2\left\|v^j_i(x)\right\|_{L^2(B_1)}^{2}+o(\epsilon)=
$$
$$
=s^{\lambda_2-\frac{1}{p-1}}\left\| \nabla' u^j(x)\right\|_{L^2(B_1(0))}^2+o(\left\| \nabla' u^j(x)\right\|_{L^2(B_1(0))}^2)\le
$$
$$
\le \frac{s^{\lambda_2-\frac{1}{p-1}}+1}{2}\left\| \nabla' u^j(x)\right\|_{L^2(B_1(0))}^2
$$
where we used (\ref{New15}) and (\ref{skrapigt}). This implies the Proposition with 
$$
\tau=\frac{s^{\lambda_2-\frac{1}{p-1}}+1}{2}<1.
$$
\qed

\begin{lem}\label{wtfaid}
 Let $u$ be a solution to the $p-$harmonic obstacle problem in a normalized coordinate system. Then there exists an $\epsilon_0>0$ such that if $0\in \Gamma_u$ and 
 $\|\nabla' u\|_{L^2(B_1(0))}\le\epsilon <\epsilon_0$ then there exists a vector $\eta$ such that
 $$
 \lim_{r\to 0}\frac{u(rx)}{r^{\frac{p}{p-1}}}= \frac{p-1}{p} (\eta\cdot x)_+^{\frac{p}{p-1}}.
 $$
 Furthermore, there exists a constant, $C$, such that
 \begin{equation}\label{SomeStars}
 \left|e_n-\eta\right|\le C \epsilon.
 \end{equation}
\end{lem}
\textsl{Proof:} By Proposition \ref{improvement} it follows that
$$
\left\| \nabla' \frac{u(sx)}{s^{\frac{p}{p-1}}} \right\|_{L^2(B_1(0))}\le \tau \epsilon.
$$
We may thus re-normalize the coordinate system (rotate the coordinate axes to a new set of basis vectors $\{e_1^1,e_2^1,...,e_n^1\}$) so that
$$
u_s=\frac{u(sx)}{s^{\frac{p}{p-1}}}
$$
satisfies
$$
\left\| \nabla' u_s\right\|_{L^2(B_1(0))}\le \tau \epsilon
$$
in the new coordinate system. 

Next we claim that $|\cos^{-1}(e_n\cdot e_n^1)|\le C\tau \epsilon$. The argument is elementary so we will only give a rough sketch. From Lemma \ref{xncontrol}
we know that $\partial_n u_s \approx (x_n)_+^{1/(p-1)}+o(1)$. We also know that $\|\nabla' u_s\|_{L^2(B_1)}\le C_s \epsilon$
since $\|\nabla' u\|_{L^2(B_1(0))}\le \epsilon$. Therefore 
\begin{equation}\label{AlgOrNot}
C\tau \epsilon \ge \|\nabla_{e_n^1} u_s\|_{L^2(B_1(0))}\ge (1-o(1))\left\|\nabla_{e_n^1}\frac{p-1}{p}(x_n)_+^{p/(p-1)} \right\|_{L^2(B_1(0))}-C_s\epsilon.
\end{equation}
A direct calculation, and using a Taylor expansion, shows that 
\begin{equation}\label{FuncCalc}
\|\nabla_{e_n^1} u_s\|_{L^2(B_1(0))}\ge c_0 (1-e_n\cdot e_n^1)
\end{equation}
for $e_n^1\approx e_n$.
Using (\ref{FuncCalc}) in (\ref{AlgOrNot}) and some simple calculus calculations implies that $|\cos^{-1}(e_n\cdot e_n^1)|\le C\tau \epsilon$,
where $C$ depend on $\tau$ and on $s$.

We may thus repeat the argument and conclude that for each $k\ge 0$ there is a coordinate system $\{ e_1^k,e_2^k,...,e_n^k\}$ such that
\begin{equation}\label{flatnessimp}
\left\| \nabla' u_{s^k}\right\|_{L^2(B_1(0))}\le \tau^k \epsilon\to 0 \textrm{ as }k\to \infty
\end{equation}
and the rotation of of $\{e_1^{k-1},...,e_n^{k-1}\}$ with respect to $\{e_1^{k},...,e_n^{k}\}$ is of order $C\tau^k \epsilon$, that is
$|\cos^{-1}(e_n^{k-1}\cdot e_n^k)|\le C\tau^k \epsilon$.

In particular,
\begin{equation}\label{IndeedSire}
|\cos^{-1}(e_n \cdot e_n^k)|\le \sum_{l=1}^{k}|\cos^{-1}(e_n^{l-1} \cdot e_n^{l})|\le C\epsilon \sum_{l=1}^k \tau^l\le C\epsilon.
\end{equation}
That means that the coordinate system converges: $e_n^k\to \eta$ for some vector $\eta$. And (\ref{flatnessimp}) implies that
$$
\gamma\cdot \nabla u_{s^k}\to 0  
$$
for any $\gamma$ orthogonal to $\eta$, $\gamma\cdot \eta=0$.

That means that $u_{s^k}\to u_0$ where $u_0$ depend only on the $\eta$ direction. That is $u_0(x)=\frac{p-1}{p} (\eta\cdot x)_+^{\frac{p}{p-1}}.$

The estimate (\ref{SomeStars}) follows directly from passing to the limit $k\to \infty$ in (\ref{IndeedSire}). \qed

\begin{cor}\label{somecor}
 Let $u$ be as in Lemma \ref{wtfaid}, maybe with a smaller $\epsilon_0$, then $\Gamma_u$ is a $C^{1,\alpha}$ graph in $B_{1/2}(0)$ for some $r_0>0$. Furthermore
 the $C^{1,\alpha}-$norm of the graph is bounded by $C\epsilon$.
\end{cor}
\textsl{Proof:} Notice that if $u$ is as in the Corollary and $x^1\in \Gamma_u\cap B_{1/2}(0)$ then
$u_{x^1,1/2}=\frac{u(x/2+x^1)}{(1/2)^{p/(p-1)}}$ will satisfy the hypotheses of Lemma \ref{wtfaid} in $B_1(0)$ with $C_n\epsilon$
in place of $\epsilon$. This implies, if $\epsilon$ is small enough, that the blow-up of $u_{x^1,1/2}$ is unique:
$$
\lim_{r\to 0}\frac{u_{x^1,1/2}(rx)}{r^{p/(p-1)}}= \frac{p-1}{p}(\eta_{x^1}\cdot x)_+^{\frac{p}{p-1}}.
$$
This defines the measure theoretic normal of $\Gamma_u$ at the point $x^1\in \Gamma_u$. Furthermore, by (\ref{SomeStars}), 
$|\eta_{x^1}-\eta_0|\le C\epsilon$. Therefore the normal of $\Gamma_u$ is well defined at every point in $B_{1/2}(0)$. 

It only remains to show that $\eta_x$ is H\"older continuous in $x$. We will show that 
$$
|\eta_{x^1}-\eta_0|\le C |x^1|^\alpha \epsilon
$$
with $\alpha=\frac{\ln(\tau)}{\ln(s)}>0$. A similar estimate for of $|\eta_x-\eta_y|$ for two arbitrary points $x,y\in \Gamma_u$ follows by a simple translation 
translating $y\in B_{r_0}(0)$ to the origin.

To that end assume that $s^{k+1}<|x^1|\le s^k$. Then $\frac{x^1}{s^k}\in \Gamma_{u_{s^k}}$ and 
$u_{s^k}$ satisfies the hypothesis of Lemma \ref{wtfaid} with $\tau^k \epsilon$ in place of $\epsilon$. This implies, again by (\ref{SomeStars}), that
$$
|\eta_{x^1}-\eta_0|\le C\tau^k \epsilon\le C |x^1|^{\frac{\ln(\tau)}{\ln(s)}}\epsilon.
$$
\qed

\begin{thm}\label{MainResult}
 Let $u$ be a solution to the $p-$harmonic obstacle problem in $B_1(0)$ then there exists an open set $\Gamma_0\subset \Gamma$
 such that $\H^{n-1}(\Gamma\setminus \Gamma_0)=0$ and for every $x^0\in \Gamma_0$ there exists an $r=r(x^0)$ such that 
 $\Gamma_u\cap B_r(x^0)$ is a $C^{1,\alpha}$ graph.
\end{thm}
\textsl{Proof:} Since $\Omega_u$ has finite perimeter, by Lemma \ref{ThisToo}, it follows that for $\H^{n-1}$-a.e. free boundary 
point $x^0\in \Gamma_u$ blow-up 
$$
\lim_{r\to 0} \frac{u(rx+x^0)}{r^{p/(p-1)}}=u^0
$$
has support in a half space which we may assume (after possibly rotation the coordinate system) to be $\R^n_+$.

By Proposition \ref{flatnessprop}  we know that for a.e. $y^0$ in $\Gamma_{u^0}$
$$
\lim_{r\to 0} \frac{u(rx+y^0)}{r^{p/(p-1)}}\to \frac{p-1}{p}(\eta\cdot x)_+^{\frac{p}{p-1}}
$$
for some vector $\eta$. That is 
$$
\frac{p-1}{p}(\eta\cdot x)_+^{\frac{p}{p-1}}\in Blo(u^0,y^0).
$$
From Lemma \ref{BlowOfBlowAreBlow} we can conclude that for $\H^{n-1}-$a.e. $x^0\in \Gamma_u$
$$
\frac{p-1}{p}(\eta\cdot x)_+^{\frac{p}{p-1}}\in Blo(u,x^0).
$$

If we normalize the coordinate system this means that $\frac{u(rx+x^0)}{r^{p/(p-1)}}$
satisfies the conditions of Corollary \ref{somecor} if $r$ is small enough. The Theorem follows. \qed

\bibliographystyle{plain}
\bibliography{p-Harmonic.bib}

\end{document}